\documentclass[10pt,twoside]{article}
\usepackage[latin1]{inputenc}
\usepackage[english,francais]{babel}
\usepackage{amsfonts}
\usepackage{amssymb}
\usepackage[all]{xy}
\usepackage{amsmath}
\usepackage{xypic}
\usepackage{texdraw}
\usepackage{amstext}
\usepackage{amsbsy}
\usepackage{amsopn}
\usepackage{amsthm}
\usepackage{amsxtra}
\usepackage{euscript}
\usepackage{graphics}
\usepackage{pstricks,pst-3d,pst-char,pst-coil,pst-eps,pst-fill,pst-grad,pst-node,pst-plot,pst-text,pst-tree}
\title{Indices d'un opérateur différentiel matriciel et applications}
\author{Kamel  BETINA\\
USTHB, Faculté de Mathématiques, BP 32, El Alia, Bab Ezzouar, Alger\\
e-mail : kamelbetina@yahoo.fr}
\date{}
\begin{document}
\maketitle
\thispagestyle{empty}
\abstract{Dans ce travail, on détermine l'indice formel et l'indice polynômial d'un opérateur différentiel linéaire  
$P = \sum\limits_{i=0}^{m} A_{i}(x)\dfrac{d^i}{dx^{i}}$ à coefficients dans $M_n(\mathbb{C}[x])$
 avec $\det A_m(x)$ non identiquement nul.\\
 Ensuite, on applique ces résultats pour donner une nouvelle démonstration d'un théorème de Bezivin-Robba équivalent au théorème de Lindemann-Wierstrass, ainsi que pour trouver des conditions suffisantes sur  l'équation différentielle matricielle de Riccati $Y'=A(x)+B(x)Y+YC(x)Y$ à coefficients dans $M_n(\mathbb{C}[x])$ pour que toute solution méromorphe soit rationnelle et d'autres conditions suffisantes pour que la solution générale soit algébrique.
 \\[0.7cm]
In this paper, one determines the formal index and the   polynômial index of a linear differential operator 
$P= \sum\limits_{i=0}^{m} A _ {i}(x)\dfrac{d^i}{dx^{i}}$ with coefficients in $M_n(\mathbb{C}[x])$ and $ \det A_m(x)\not\equiv 0$ .\\
Then, one applies these results to give a new proof of a   Bezivin-Robba theorem equivalent to the  Lindemann-Wierstrass theorem, as to find  sufficient conditions on the Riccati  matrix differential equation    $Y'=A(x)+B(x)Y+YC(x)Y$ with coefficients in $M_n(\mathbb{C}[x])$ so that  any meromorphic solution is rational and other sufficient conditions  so that the general solution is algebraic.}\\[0.25cm]
\textit{AMS classification subject} : 34M, 11J\\[0.25cm]
\section{Notations}
\hspace{-.5cm}$\mathbb{C}[x]$=l'anneau des polynômes en $x$ à coefficients dans $\mathbb{C}$,\\
$\mathbb{C}[x]_a=\{\dfrac{b}{a^n}\, ;\,  b \in \mathbb{C}[x],n \in \mathbb{N}  \}$   = le localisé en $a\in \mathbb{C}[x]$ de l'anneau $\mathbb{C}[x]$,\\
$\mathbb{C}[[x]]$=l'anneau des séries formelles en $x$ à coefficients dans $\mathbb{C}$,\\
$\mathcal{M}(\mathbb{C})$ = le corps des fonctions méromorphes sur $\mathbb{C}$,\\
Pour $Y=\{a_1,...,a_n\} \subset \mathbb{C}$, on pose $\mathcal{M}_Y(\mathbb{C})=\{f\in \mathcal{M}(\mathbb{C})\,; f\, \text{est analytique dans} \, \mathbb{C}-Y\}$.\\ 
Si $g\in \mathbb{C}[x]$, on note $ord_{\alpha}g$ l'ordre du zéro de $g$ au point $\alpha$ et $d°g$ le degré du polynôme $g$ ; si $A=(a_{ij})$ est une matrice à coefficients dans $\mathbb{C}[x]$, on note $d°A$ le nombre $\sup\limits_{ij}d°a_{ij}$.\\
Si $u:E\rightarrow E$ est un endomorphisme d'espaces vectoriels, on note $ker(u,E)$ et $Coker(u,E)$ respectivement son noyau et son conoyau.\\
On note $v(f)$ la valuation de $f=\sum\limits_{i\geq i_0}f_ix^i \in \mathbb{C} [[x]][x^{-1}]$ :$$v(f)=\begin{cases}
\inf \{i;f_i\ne 0 \}, & \text{si $f\ne 0$;}\\
+\infty & \text{si $f=0$ }
\end{cases}$$
Si $A=(a_{ij}) \in End(\mathbb{C}[x]^N)$, on pose $\nu (A)=inf\{v(a_{ij}); 1\leq i,j \leq N \}$. \\[0.2cm]
 On considère l'opérateur différentiel matriciel  $$P = \sum\limits_{i=0}^{m} A_{i}(x){d^i \over dx^{i}}$$
  avec 
   $$ A_{i}(x)=A_{i}^0+A_{i}^{1}x+
  ... +A_{i}^{d_i}x^{d_i}\in M_{n}(\mathbb{C}[x]^N) \hspace{0.15cm}  \text{et} \hspace{0.15cm}\det
  A_{m}(x)\not\equiv 0 .$$ \\
  Si $A_i(x)\equiv 0$, on pose $d_i=-\infty$.
  \newline\\
  Soient 
  $$n_i=\inf \{j;A_i^j\ne 0 \}=\nu (A_i),$$ 
  $$n=\sup\limits_ {0\leq i \leq
  m}(i-n_i)=\sup\limits_ {0\leq i \leq
  m}(i-\nu (A_i)),$$  
  $$ J=\{j ; j-n_j =n\},$$
  $$ L(k)=\sum\limits_{j \in J}k(k-1) ... (k-j+1)A_{j}^{n_j}
 ,$$
  $$ l(k)=\det L(k).$$ 
  Remarquons que si $m \in J$, alors $l(k)$ est un polynôme en $k$ de
  degré $mN$. \newline\\
  Soient 
  $$n'=\inf\limits_{0\leq i \leq m}(i-d_i),$$ 
  $$J'=\{j ; j-d_j=n'\},$$
  $$D(k)=\sum\limits_{j \in J'}k(k-1) ... (k-j+1)A_{j}^{d_j},$$  
  et 
  $$d(k)=\det D(k).$$
  \section{Indice formel et indice polynômial}
\hspace{-.5cm}\underline{\bf Définition}  :\\  \textit{Le polynôme } $l(k)$  \textit{est appelé polynôme
  indiciel de l'opérateur matriciel }$P$ \textit{et l'\'{e}quation } $l(k)=0$
  \textit{est appel\'{e}e l'\'{e}quation indicielle de l'\'{e}quation
  diff\'{e}rentielle matricielle }$Pu=0$.\\[0.5cm]
\underline{\bf Définition} :\\ \textit{Soit $a\in \mathbb{C}$. L'op\'{e}rateur différentiel
matriciel} $P$ \textit{est dit
 régulier  en } $a$ \textit{si toute solution de l'équation $Pu=0$ est analytique}, \textit{il est dit singulier régulier en $a$ si toute d\'{e}termination $f$ de toute
  solution de l'\'{e}quation} $Pu=0$ \textit{est à croissance
  mod\'{e}r\'{e}e au voisinage de} $a$ \,\textit{i.e si pour tout secteur $S$ de sommet $a$ il existe $\beta \in \mathbb{N}$ tel que} $\lim\limits_{x \to a \atop x\in S} (x-a)^{\beta}f=0$, \textit{ $P$ est dit singulier irrégulier en} $a$ \textit{s'il n'est ni régulier ni singulier régulier en $a$}. \newline\\
  \underline{\bf Définition}:\\ \textit{Soient $E$,\,$F$ deux espaces vectoriels sur un corps $K$ et $u:E \rightarrow F$ un homomorphisme. On dit que l'opérateur $u$ est à indice si $\dim Keru$ et $\dim Cokeru$ sont finis et on apelle alors indice de} $u:E\rightarrow F$ \textit{le nombre entier $\chi (u;E,F)=\dim Keru-\dim Cokeru$. Si $E=F$, on note $\chi (u;E)$ l'indice $\chi (u;E,F)$ ( s'il existe )}.\\[0.2cm]
  D'après [5], l'opérateur $P:\mathbb{C}^N\{x\}^N \rightarrow \mathbb{C}\{x\}^N$ est à indice et on a : $$ \chi(P;\mathbb{C}\{x\}^N)=mN-v(detA_m) \hspace{3cm} (1)$$
  D'après [1], l'opérateur $P:\mathbb{C}[[x]]^N\rightarrow \mathbb{C}[[x]]^N$ est à indice et les quatre conditions suivantes sont équivalentes \\
  (i)\; $P$ est singulier r\'{e}gulier en $0$ ;\\
  (ii)\; $\chi(P,\Bbb{C}[[x]]^N)=\chi(P,\Bbb{C}\{x\}^N)$ ;\\
 (iii)\; l'opérateur $$ P :{\mathbb{C}[[x]]^N \over \Bbb{C}\{x\}^N}\rightarrow {\mathbb{C}[[x]]^N \over
  \mathbb{C}\{x\}^N}$$
  est injectif ;\\
  (iv)\; l'opérateur différentiel d'ordre $1$ suivant est singulier régulier à l'origine\newline\\
  $$ D={d \over dx} \; - \left (
  \begin{array}{ccccc}
  0&I&0&\cdots &0 \\
  0&0&I& \cdots &0\\
  \vdots & \vdots & \ddots &\ddots & \vdots \\
  0&0&0&\ddots &I\\
  B_0& B_1& B_2& \cdots & B_{m-1}
\end{array}
  \right )
  $$
  où $B_i=A_m^{-1}A_i \; (0\leq i \leq m-1)$ \, et $I = $\,  la
  matrice identité d'ordre $N$.\\[0.3cm]
  Dans la suite, on travaillera au voisinage du point $a=0$. 
\\[0.2cm]
  \underline{\bf Théorème $1$} :\\[0.2cm]
$(i)$) \textit{Si} $l(k)\not\equiv 0$, \textit{alors} $\chi (P;\mathbb{C}[[x]]^N)=nN$ ,\\
$(ii)$) \textit{Si} $m \in J$, \textit{alors 
$P$ est singulier régulier en $0$ si, et seulement si}, $ v(\det (A_m))=\nu (A_m)N$.
\\[0.3cm]
\underline{Preuve}:\\
$(i)$ Puisque $\sup\limits_i (i-v(A_i))=n$, on a :\\
$P(f_0x^k)=L(k)f_0x^{k-n}+L^1(k)f_0x^{k-n+1}+ ...$\\[0.2cm]
$P(f_1x^{k+1})= L(k+1)f_1x^{k-n+1}+ ...$\\
\vdots
\\[0.3cm]
De là, on déduit par un calcul simple, un isomorphisme 
 $$P:x^k\mathbb{C}[[x]]^N \xrightarrow{\thicksim}x^{k-n}\mathbb{C}[[x]]^N $$
 et du diagramme commutatif :\\
$$\xymatrix{
0\ar[r]&x^k\mathbb{C}[[x]]^N \ar[d]^{P} \ar[r]&\mathbb{C}[[x]]^N \ar[d]^
P\ar[r]&\mathbb{C}[[x]]^N/x^k\mathbb{C}[[x]]^N \ar[r] \ar[d]^P&0 \\
0\ar[r]& x^{k-n}\mathbb{C}[[x]]^N\ar[r]& \mathbb{C}[[x]]^N\ar[r]& \mathbb{C}[[x]]^N/x^{k-n}\mathbb{C}[[x]]^N\ar[r]&0
}
$$
on déduit, à l'aide du lemme du serpent, que 
 $$\chi(P;\mathbb{C}[[x]]^N)=kN-(k-n)N=nN$$
 d'où la première partie du théorème.\\
 $(ii)$ Si $m\in J$, $l(k)$ est un polynôme en $k$ de degré $N\sup J=Nn$, donc $l(k)\not\equiv 0$,  d'où: 
\begin{eqnarray*}
 \chi(P;\mathbb{C}[[x]]^N)&=&nN \hspace{1.5cm} (\text{d'après $(i)$ })\\
  &=&\sup\limits_i\{i-\nu (A_i)\}N\\
   &=&(m-\nu (A_m))N  \hspace{1.5cm} (\text{car $m\in J$}).\\
\end{eqnarray*}
Si 
$$v(\det A_m)=\nu (A_m)N,$$ 
alors 
$$(m-\nu (A_m))N=mN-v(\det (A_m)),$$
 i.e
 $$\chi(P;\mathbb{C}[[x]]^N=\chi(P;\mathbb{C}\{x\}^N,$$
  c'est à dire que $P$ est singulier régulier en $0$.\\
Réciproquement, si $P$ est singulier régulier en $0$, alors : $$\chi(P;\mathbb{C}[[x]]^N=\chi(P;\mathbb{C}\{x\}^N),$$
 mais comme $l(k) \not\equiv 0$, on a: 
 $$\chi(P;\mathbb{C}[[x]]^N=nN\hspace{1cm}\text{(d'après  $(i)$)},$$ 
et
$$\chi(P;\mathbb{C}\{x\}^N)=mN-v(det(A_m))\hspace{1cm}  \text{(d'après [$5$])}$$
donc 
$$nN=(m-\nu (A_m))N \hspace{1cm} \text{(car $m\in J$)},$$ \\
d'où l'égalité\\
$$v(det(A_m))=\nu (A_m)N.$$
\underline{\bf Théorème $2$} :\\
$1)$ \textit{Si}\quad  $m \notin J$ \textit{et} $l(k)\not\equiv 0$, \textit{alors $P$ est singulier irrégulier en $0$}.\\
$2$ \textit{Si $m \notin J$ et $l(k)\equiv 0$, $P$ peut être
singulier régulier ou singulier irrégulier}.\\[0.2cm]
\underline{Preuve}:\\
1) Si $m\notin J$, alors  $m-v(A_m)<n$, i.e$$(m-v(A_m))N<nN \hspace{4cm} (2)$$
Si $l(k)\not\equiv 0$, alors $\chi(P;\mathbb{C}[[x]]^N)=Nn$ (Théorème $1$ ).\\
Donc, si $m\notin J$ et $l(k)\not\equiv 0$, on a : 
\begin{eqnarray*}
 \chi(P;\mathbb{C}[[x]]^N)-\chi(P;\mathbb{C}\{x\}^N)&=&Nn-(m-v(A_m))N>0  \hspace{2cm} (\text{d'après $(2)$})
\end{eqnarray*}
i.e $P$ est singulier irrégulier.
\\[0.2cm]
2) En effet : \\
$a$)\; Si $$P_1=x^2I_2{d \over dx} +  \begin{pmatrix}
0&0\\
1&-1 \end{pmatrix}  \hspace{1.5cm} (\text{où $I_2$ est la matrice identité d'ordre $2$})$$
on a :\\
$l(k)\equiv 0$ ,\\
$m=1$ ,\\
et\\
$n=\sup\{0-0,1-2\}=0$ ,\\
donc $m \notin J$.\\
Le vecteur
$u=\left ( 0 \; , \; \sum\limits_{n\geq 0}n!x^{n+1} \right )^t$
  \, vérifie l'équation $$P_{1}u = \left (
  \begin{array}{c}
0 \\
-x
\end{array}
\right ) \, ,$$
 donc l'opérateur $$ P_1 :{\mathbb{C}[[x]]^N \over \Bbb{C}\{x\}^N}\rightarrow {\mathbb{C}[[x]]^N \over
  \mathbb{C}\{x\}^N}$$
  n'est pas injectif,
  i.e $P_1$ est singulier irrégulier en $0$.\\
 $b$)\; Si $$P_2=x^3I_2{d \over dx} +  \begin{pmatrix}
0&0\\
1&0 \end{pmatrix}  $$ 
on a :\\
$l(k)\equiv 0$ ,\\
$m=1$,\\
et\\
$n=sup\{0-0,1-3\}=0$ ,\\
donc $m \notin J $.\\[0.2cm]
Si $v=\left ( \begin{array}{c}
v_1\\v_2
\end{array}
\right )$ $\in \mathbb{C}[[x]]^2$ vérifie l'équation $$P_2v=g=\left ( \begin{array}{c}
g_1\\g_2
\end{array}
\right ) \in \mathbb{C}\{x\}^2 \, ,$$  alors on a le système différentiel 
$$
\left \lbrace
\begin{array}{l}
x^3v'_1=g_1\\
v_1+x^3v'_2=g_2
\end{array}
\right.
$$
d'où $v_1$ ,$v_2$ $\in \mathbb{C}\{x\}^2$ et donc l'application 
$$
P_2 :{\mathbb{C}[[x]]^N \over \Bbb{C}\{x\}^N}\rightarrow {\mathbb{C}[[x]]^N \over
  \mathbb{C}\{x\}^N}
$$
est injective, i.e $P_2$ est singulier régulier.\\[0.2cm]
\underline{\bf Théorème $3$} :\\[0.2cm]
\textit{Si $l(k)\not\equiv 0$ et $d(k)\not\equiv 0$, alors }
$$
\chi(P;\mathbb{C}[x]^N)=Nn'=N\inf\limits_i(i-d°A_i) .
$$
\underline{Preuve}:\\
Pour $k\in \mathbb{N}$ , on a :\\
$\dfrac{d^i}{dx^i}(x^k)=k(k-1) ... (k-i+1)x^{k-i}$\\
et donc, pour $\lambda \in \mathbb{C}^N$,\\
\begin{eqnarray*}
P(\lambda x^k)&=& \sum\limits_{i=0}^mA_i\dfrac{d^i}{dx^i}(\lambda x^k)\\
 &=&L(k)\lambda x^{k-n}+L_1(k)\lambda x^{k-n+1}+ ... +L_{n-n'-1}(k)\lambda x^{k-n'-1}+ D(k)\lambda x^{k-n'}
\end{eqnarray*}
Ainsi pour un polynôme \\
$f=f_0x^k+ ... +f_lx^{k+l} \in \mathbb{C}[x]^N$\\
on obtient 
\begin{eqnarray*}
Pf&=&L(k)f_0x^{k-n}+L_1(k)f_0x^{k-n+1}+ ...+L_{n-n'-1}(k)f_0x^{k-n'-1}\\
 & &+ ... +L(k+l)f_lx^{k-n+l}+ .. +D(k+l)f_lx^{k-n'+l}
 \end{eqnarray*}
 \\Désignons par $M(k)$ l'espace des polynômes de $\mathbb{C}[x]$ de valuation supérieure ou égale à $k$.\\[0.32cm]
 \underline{\bf lemme $1$}:\\[0.2cm]
 Pour $k$ assez grand, on a :
$$\dim Coker(P:M(k)^N\rightarrow M(k-n)^N)=(n-n')N$$\\[0.2cm]
\underline{preuve}:\\
Pour $k$ assez grand, $l(k)=\det L(k)$ et $d(k)=\det D(k)$ ne s'annulent pas.\\
 Si $n=n'$, l'application $P:M(k)^N\rightarrow M(k-n)^N$ est manifestement surjective.\\
 Si $n\ne n'$, soit $\psi_i=(\delta_i^1,...,\delta_i^N)^t$ où $\delta _i^j$ est le symbole de Kronecker ($1\leq i\leq N$). Les $(n-n')N$ polynômes $x^{k-n}\psi_i$, ... ,$x^{k-n'-1}\psi_i$ forment une base d'un espace supplémentaire de $P(M(k)^N)$ dans $M(k-n)^N$. En effet, ils n'appartiennent pas à $P(M(k)^N)$ d'une part, d'autre part, tout $g\in M(k-n)^N$ s'écrit :\\
$g=g_0x^{k-n}+...+g_{n-n'+l}x^{k-n'+l} \quad , g_i\in \mathbb{C} ^N\quad (0\leq i\leq n-n'+1)$\\
et il existe deux uniques polynômes à coefficients dans $ \mathbb{C}^N$\\
$f=f_0x^k+...+f_lx^{k+l}$ et $v=v_0x^{k-n}+...+v_{n-n'-1}x^{k-n'-1}$\\
tels que $Pf=g-v$.\\
Le polynôme $f$ est déterminé (pour $k$ assez grand) par les $l+1$ équations :$$
\begin{array}{l}
D(k+l)f_l=g_{n-n'+l}\\
D(k+l-1)f_{l-1}+L_{n-n'-1}(k+l)f_l=g_{n-n'+l-1}\\
\vdots \\
D(k+l-(n-n'+1))f_{l-(n-n'+1)}+...+L(k+l)f_l=g_{l-1}\\
\vdots \\
L(k+n-n'))f_{n-n'}+...+D(k)f_0=g_{n-n'}
\end{array}$$
et $v$ par les $(n-n')$ équations 
$$
\begin{array}{l}
L(k)f_0=g_0-v_0\\
L_1(k)f_0+L(k+1)f_1=g_1-v_1\\
\vdots \\
L_{n-n'-1}(k)f_0+...+L(k+n-n'-1)f_{n-n'-1}=g_{n-n'-1}-v_{n-n'-1}
\end{array}
$$
d'où le lemme $1$.\\
Pour $k$ assez grand, l'application $P:M(k)\rightarrow M(k-n)$ est injective et donc du diagramme commutatif :\\
$$\xymatrix{
0\ar[r]&M(k)^N \ar[d]^{P} \ar[r]&\mathbb{C}[x]^N \ar[d]^
P\ar[r]&\mathbb{C}[x]^N/M(k) \ar[r] \ar[d]^P&0 \\
0\ar[r]& M(k-n)\ar[r]& \mathbb{C}[x]^N\ar[r]& \mathbb{C}[x]^N/M(k-n)\ar[r]&0
}
$$
on déduit que :\\
$\chi(P,\mathbb{C}[x])=\chi(P,M(k),M(k-n))+kN-(k-n)N=(n'-n)N+nN=n'N$.\\[0.2cm]
\underline{\bf Théorème $4$}:\\[0.2cm]
On suppose que $l(k)\not\equiv 0$ et $d(k)\not\equiv 0$. Si $Nn'=mN-d°(det(A_m(x)))$, alors 
$$
(Pu\in \mathbb{C}(x)^N, u\in \mathcal{M}(\mathbb{C})^N)\Rightarrow u\in \mathbb{C}(x)^N
$$
\underline{Preuve} :\\
Il suffit de prouver que l'application 
$$
P :{\mathcal{M}(\mathbb{C})^N \over \mathbb{C}(x)^N} \rightarrow  {\mathcal{M}(\mathbb{C})^N \over \mathbb{C}(x)^N}
$$
est injective.\\
Soit $u \in \mathcal{M} (\mathbb{C})^N $ vérifiant l'équation 
$$
Pu=g \in \mathbb{C}(x)^N \quad ,
$$
$g$ s'écrit sous la forme 
$$
g=\frac{1}{s}g_1  \hspace{1.5cm} \text{avec $s\in \mathbb{C}[x]$ et $g_1 \in \mathbb{C}[x]^N$} .
$$
En dehors de l'ensemble $Y=\{z\in \mathbb{C}; s(z)detA_m(z)=0\}$, le vecteur $u$ est analytique. On a donc
$$
u \in \mathcal{M}_Y(\mathbb{C})^N
$$
et, en posant $t=t(z)=s(z)detA_m(z)$, le problème revient donc à montrer le lemme suivant :\\
\underline{\bf lemme $2$}:\\
\textit{Sous les hypothèses du} Théorème $4$, \textit{l'application} 
$$
P :{\mathcal{M}_Y(\mathbb{C})^N \over \mathbb{C}[x]_t^N} \rightarrow  {\mathcal{M}_Y(\mathbb{C})^N \over \mathbb{C}[x]_t^N}
$$
\textit{est injective}.\\[0.2cm]
\underline{Preuve}:\\
D'après [5], on a :
$$
\chi(P;\theta (\mathbb{C})^N)=mN(\dim H^0(\mathbb{C},\mathbb{C})-\dim H^1(\mathbb{C},\mathbb{C}))-\sum\limits_{x\in \mathbb{C}}ord_x(\det (A_m))\\
$$
or\\
$H^0(\mathbb{C},\mathbb{C})\simeq \mathbb{C}$ ,\\
$H^1(\mathbb{C},\mathbb{C})=0$ ,\\
et\\
$\sum\limits_{x\in \mathbb{C}}ord_x(\det A_m(x))= d°(\det A_m(x))$,
\\
d'où \\[0.2cm]
$
\chi(P;\theta (\mathbb{C})^N)=mN-d°det(A_m(x)).
$\\[0.2cm]
$\mathbb{C}[x]$ étant dense dans $\theta (\mathbb{C})$, donc l'application 
$$
P :{\theta (\mathbb{C})^N \over \mathbb{C}[x]^N} \rightarrow  {\theta (\mathbb{C})^N \over \mathbb{C}[x]^N}
$$
est surjective; si on note $A_{x_i}$ l'espace vectoriel engendré par la famille $\left(\dfrac{1}{(x-x_i)^j}\right)_{j\geq 1}$ \quad , on a :\\[0.3cm]
$$
\mathcal{M}_Y(\mathbb{C})=\theta (\mathbb{C})\underset{x_i \in Y}{\oplus} A_{x_i}\quad ,
$$
et
$$
\mathbb{C}[x]_t=\mathbb{C}[x]\underset{x_i \in Y}{\oplus} A_{x_i} \quad ,
$$
d'où 
$$
{\mathcal{M}_Y(\mathbb{C})^N \over \mathbb{C}[x]_t^N}\simeq {\theta (\mathbb{C})^N \over \mathbb{C}[x]^N}  \hspace{1cm} (3)
\quad ,$$
par suite 
$$
Coker\left( P,{\mathcal{M}_Y(\mathbb{C})^N \over \mathbb{C}[x]_t^N}\right) =(0)
$$
et donc
\begin{eqnarray*}
\dim ker\left( P,{\mathcal{M}_Y(\mathbb{C})^N \over \mathbb{C}[x]_t^N}\right)&=&\chi \left( P;{\mathcal{M}_Y(\mathbb{C})^N \over \mathbb{C}[x]_t^N}\right)\\
 &=&\chi \left( P;{\theta(\mathbb{C})^N \over \mathbb{C}[x]^N}\right) \hspace{1.5cm} (\text{d'après}\, (3))\\
 & =& \chi(P;\theta(\mathbb{C}^N))-\chi(P;\mathbb{C}[x]^N) \\
  &=&mN-d°(detA_m(x))-\inf\limits_i (i-d°A_i)N\\
   &=& mN-d°(detA_m(x))-n'N\\
   &=&0  \hspace{2cm}  (\text{d'après l'hypothèse du Théorème $4$} )
\end{eqnarray*}
d'où le lemme $2$.
\section{Applications}
\subsection{Nouvelle démonstration du théorème de Lindemann-Wierstrass}
\hspace{-0.5cm}\underline{\bf Théorème} (Lindemann-Weierstrass):\\[0.3cm]
\textit{Soient $a_1,a_2,\, ...\, ,a_n$ des nombres algébriques deux à deux disjoints et $\lambda_1,\lambda_2,\, ...\, , \lambda_n$ des nombres algébriques non tous nuls, alors}  $\sum\limits_{i=1}^{s}\lambda_i e^{a_i}\ne 0$.\\[0.3cm]
Dans[4], Bezivin et Robba ont démontré que ce dernier résultat est équivalent au théorème suivant :\\[0.2cm]
\underline{\bf Théorème $5$}:\\[0.3cm]
\textit{Soit $L=x^2\dfrac{d}{dx}+x-1$. Si $u\in \mathbb{C}\{x\}$ vérifie $Lu \in \mathbb{C}(x)$, alors} $u \in \mathbb{C}(x)$.\\[0.3cm]
Ensuite, ils ont prouvé ce dernier théorème avec des méthodes $p$-adiques.\\
Maintenant, on va donner une nouvelle démonstration du Théorème $5$ en utilisant le théorème $3$ ainsi que la théorie classique des équations différentielles dans le champ complexe.\\
En posant $P=L$, on a :\\
$N=m=1 , \\
d(k)=k+1, \\
l(k)\not \equiv 0,\\
n'=inf(0-1,1-2)=-1,$\\
et\\
$d°detx^2=2$.\\
Ainsi \, $d(k)\not\equiv 0$\, , $l(k)\not\equiv 0$ , $Nn'=-1$ \, et $mN-d°detx^2=1-2=-1$, donc, d'après le \mbox{ Théorème $4$}, si on montre que $u \in \mathcal{M}(\mathbb{C})$, on pourra conclure que $u\in \mathbb{C}(x)$, ce qui terminera la preuve du Théorème $5$.\\
Le lemme suivant nous permettra donc de conclure.\\[0.2cm]
\underline{\bf lemme $3$}:\\[0.2cm]
\textit{Sous les hypothèses du Théorème $5$, on a} 
$
u \in \mathcal{M}(\mathbb{C})
$.\\[0.2cm]
\underline{Preuve}:\\[0.1cm]
Soit $g=Lu\in \mathbb{C}(x)$. Puisque $u\in \mathbb{C}\{x\}$ , $g\in \mathbb{C}(x)$ et donc $g=\dfrac{a}{b}$ avec $a,b \in \mathbb{C}[x]$ et $b(0)\neq 0$.\\
Soit $Z=\{x\in \mathbb{C}; \, b(x)=0\}$. Au voisinage de $0$, $u$ est analytique et au voisinage de tout autre point de $\mathbb{C}-(Z\cup \{0\})$, la solution générale de l'équation différentielle 
$$
Lu=g \hspace{3cm} (F)
$$
est analytique. Comme en tout point de $Z$ l'opérateur $L$ est régulier, donc la solution de ($F$) est à croissance modérée, il nous reste donc à vérifier qu'au voisinage de tout point de $Z$, la solution générale de ($F$) est uniforme pour conclure que $u\in \mathcal{M}(\mathbb{C})$.\\
Au voisinage d'un point $x_0\in Z$, l'équation ($F$) se transforme en 
$$\dfrac{du}{dx}+\alpha u=\dfrac{c}{(x-x_0)^s}=g$$
avec $\alpha$ , $c \in \mathbb{C}\{x-x_0\}$, $s\in \mathbb{N}^*$.\\
Posons :\\[0.2cm]
$$y=x-x_0,\hspace{0.5cm} M=\dfrac{d}{dy}+\alpha ,$$
$$u=u_0y^k+u_1y^{k+1}+ \, ... \hspace{0.5cm} \text{($k<0$)},$$
$$c=c_0+c_1y+c_2y^2+ \, ....$$ 
On a :
$$Mu=g\, \Leftrightarrow \, y^sMu=c ,$$
donc 
$$\left(c\dfrac{d}{dy}-c'\right)\left(y^s\dfrac{du}{dy}+y^s\alpha u\right)=0  \hspace{1.5cm}\text{( où $c'=\dfrac{dc}{dy}$ )}.$$
Si on pose $$P=\left(c\dfrac{d}{dy}-c'\right)(y^sM)=\left(c\dfrac{d}{dy}-c'\right)\left(y^s\dfrac{d}{dy}+y^s\alpha \right),$$ 
on a 
$$\begin{array}{rcl}
P(y^k)&=&\left(c\dfrac{d}{dy}-c'\right)(ky^{k+s-1}+\, ... )\\[0.4cm]
 &=&c_0k(k+s-1)y^{k+s-2}-c_1ky^{k+s-1}\, ...
\end{array}$$
donc les racines du polynôme indiciel de $P$ sont $k_1=0$ et $k_2=1-s$.\\
Comme $k_1-k_2=s-1 \in \mathbb{N}$, une base de $kerP$ est constituée par les solutions
$$
w_1=y^{k_1}\phi_1(y)=\phi_1(y) \quad \text{et} \quad w_2=aw_1logy+y^{k_2}\phi_2(y)
$$
où $a$ est une constante complexe, $\phi_1(y) \, \text{et} \, \phi_2(y) \in \mathbb{C}\{y\}$.\\
Toute solution $u$ de ($F$) vérifie l'équation $Pu=0$, donc $u$ est de la forme 
\begin{eqnarray*}
u&=&\gamma w_1(y)+\delta w_2(y) \hspace{2cm} (\gamma \, ,\delta \in \mathbb{C})\\
 &=&\gamma \phi_1(y)+\delta (aw_1logy+y^{k_2}\phi_2(y))
\end{eqnarray*}
Pour montrer que $u$ est uniforme au voisinage de $y=0$, il suffit de montrer que $a=0$.\\
Supposons que $a$ soit différent de $0$. On a :\\[0.2cm]
$$Pu=\left(c\dfrac{d}{dy}-c'\right)(y^sM)u=0,$$
donc
$$y^sMu=\mu c \hspace{0.3cm} \text{où $\mu \in$} \mathbb{C},$$
et
$$y^sM(a\phi_1(y)logy)=\dfrac{1}{\delta}\left( \mu c -y^sM(\gamma \phi_1(y)+\delta y^{k_2}\phi_2(y))\right)=h \in \mathbb{C}\{y\}[y^{-1}].$$
Or \\
$$y^sM(a\phi_1(y)logy)=y^s\left(a\phi '_1(y)logy+\dfrac{a\phi_1(y)}{y}+a\alpha \phi_1(y)logy \right),$$
d'où 
$$logy=\dfrac{h-ay^{s-1}\phi_1(y)}{a(\phi'_1(y)+\alpha \phi_1(y))} \in \mathbb{C}\{y\}[y^{-1}] \hspace{1.5cm} \text{(contradiction)},$$
et donc, nécessairement, $a=0$.
\\\subsection{Equation différentielle matricielle de Riccati}
On consid\`ere l'\'equation
diff\'erentielle matricielle de Riccati : $$
Y'(x)=A(x)+B(x)Y(x)+Y(x)C(x)Y(x) \qquad (R_1)$$ où $A(x),\, B(x), \,
C(x) \in M_{n}(\mathbb{C} [x]^N)$\,, $\det C(x)\not\equiv 0$,  et
$Y(x)$ est la
matrice inconnue.\\
Le changement de fonction inconnue
  $$Y(x)=-C(x)^{-1}W'(x)W(x)^{-1}$$ nous donne l'équation
  différentielle linéaire matricielle
  $$W''(x)-\{C(x)B(x)C(x)^{-1}+C'(x)C(x)^{-1}\}W'(x)+C(x)A(x)W(x)=0
  \qquad (R_2) $$
Si on pose $C^{-1}(x)=\dfrac{1}{\det C(x)} C_0(x)$ et 
$\det C(x)=c_1(x)\in \mathbb{C}[x]$, alors l'équation $(R_2)$ se transforme en \\
$$A_2(x)W'' +A_1(x)W'+ A_0(x)W=0 \hspace{1.5cm} (R_3)$$
où 
$$A_2(x)=c_1(x)I_N ,$$ 
$$A_1(x)=-\left \{C(x)B(x)C_0(x)+C'(x)C_0(x)\right \} ,$$
et 
$$A_0(x)=c_1(x)C(x)A(x).$$
On remarquera que les points singuliers de $(R_1)$ ( i.e les points où la solution générale de $(R_1)$ n'est pas holomorphe) sont les points de l'ensemble $Z=\{x\in \mathbb{C}\, ; c_1(x)=0\}$.\\[0.2cm]
\underline{\bf Théorème $6$ }:\\[0.3cm]
$1$)\textit{Si} $2-v(c_1(x))\geq \sup (1-\nu (A_1(x)), -\nu (A_0(x))$, \textit{alors}
 \textit{toute série formelle \\$Y=\sum\limits_{n\geq 0}Y_nx^n$  solution de l'équation $(R_1)$ est convergente}.\\
$2$) \textit{Si en tout point $z\in Z$, l'opérateur $P$ est singulier régulier, alors toute solution uniforme de $(R_1)$ est méromorphe}.\\[0.2cm]
\underline{Preuve}:\\[0.1cm]
$1)$ En considèrant l'opérateur différentiel 
$$P=A_2(x)\dfrac{d^2}{dx^2}+A_1(x)\dfrac{d}{dx}+A_0(x)\quad ,$$
on obtient :\\
$m=2\in J$,\,  puisque $2-v(c_1(x))\geq \sup (1-\nu (A_1(x)), -\nu (A_0(x))$ .\\
Mais comme \\
$v(\det (A_2(x)))=v(\det (c_1(x)I_N))=v(c_1^N(x))=Nv(c_1(x))=N\nu(c_1(x)I_N)=N\nu (A_2(x))$,\\
on déduit que l'opérateur $P$ est singulier régulier (Théorème $1$), donc toute solution série formelle  de l'équation $Pu=0$ est convergente.\\
Or $W$ vérifie aussi l'équation différentielle
$$W'=DW \hspace{0.4cm}\text{avec $D=-CY=D_0+D_1x+D_2x^2+\hdots \in End(\mathbb{C}[[x]]^n) $},$$
donc en posant $W=W_0+W_1x+W_2x^2+\hdots $, on obtient par identification les relations  
$$W_1=D_0W_0, $$
$$ 2W_2=D_0W_1+D_1W_0  $$
$$\vdots$$
$$nW_n=H_n(W_0,...,W_{n-1},D_0,...,D_{n-1}\hspace{0.5cm}$$
$$\vdots$$
où $H_n$ est un polynôme à $2n$ variables.\\
Ces relations nous permettent, en choisissant arbitrairement $W_0$, de calculer de proche en proche $W_1$, $W_2$ etc ..., ce qui signifie que $W$ est une série formelle solution de l'équation $Pu=0$, donc $W$ est convergente, 
par suite la série formelle $Y=-C^{-1}(x)W'W^{-1}$, qui est solution de l'équation ($R_1$), est convergente.\\[0.15cm]
$2)$ On suppose que $P$ est singulier régulier en tout point $z\in Z$, et soit $Y$ une solution uniforme  de l'équation ($R_1$). En posant $Y=-C^{-1}(x)W'W^{-1}$, $W$ est une solution de l'équation ($R_2$), donc $W$  est à croissance modérée au voisinage de chaque point singulier, par conséquent $Y$ est méromorphe,  d'où la deuxième partie du Théorème $6$.\\[0.2cm]
\underline{\bf Théorème $7$}:\\[0.1cm]
\textit{On considère le polynôme $d(k)$ associé à l'opérateur différentiel  $P$. On suppose que :\\
$1)$\,$d(k)\not\equiv 0$, \\
$2)$\, $2-d°c_1(x)\leq \inf(1-d°A_1(x), -d°A_0(x))$, \\ 
$3)$\, en tout point singulier $z\in Z$, les racines du polynôme caractéristique de $P$ sont entières, distinctes et que si une racine $r'$ est égale à $r-d$, où $r$ est une autre racine et $d$ est un entier $>0$, le calcul des coefficients de la solution (de l'équation $Pu=0$)
$$ x^{r'}+g_{r'+1}x^{r'+1}+ ...+g_rx^r+...$$
 est possible jusqu'à $g_r$.\\  
 Alors toute solution  de l'équation $(R_1)$ est rationnelle}.\\[0.2cm]
\underline{Preuve} :\\
En écrivant explicitement les coefficients de l'opérateur différentiel $$P=A_2(x)\dfrac{d^2}{dx^2}+A_1(x)\dfrac{d}{dx}+A_0(x)\quad ,$$
on a :\\
$A_2(x)=c_1(x)I_N=c_1^{n_2}I_Nx^{n_2}+\, ...\quad (\,\text{avec} \, c_1^{n_2}\ne 0 \,) \, ,$\\[0.1cm]
$A_i(x)=A_i^{n_i}x^{n_i}+
\, ... \quad  (i=0,1) \, ,$\\[0.1cm]
d'où\\
$L(k)=k(k-1)c_1^{n_2}I_N + \text{termes en $k$ de degré   $> 2$}$,\\
donc le polynôme \\
$l(k)=\det L(k)=k^{2N}(c_1^{n_2})^2+ \text{termes en $k$ de degré  $> 2$}$\\
n'est pas identiquement nul.\\
D'autre part, $d(k)\not\equiv 0$ par hypothèse et la relation \\
$2-d°c_1(x)\leq \inf(1-d°A_1(x), -d°A_0(x))$\\
 signifie que \\
$(2-d°c_1(x))N=\inf (2-d°A_2(x),1-d°A_1(x),-d°A_0(x))N=n'N$,\\
donc, d'après le Théorème $4$, toute solution méromorphe $W$ de l'équation $Pu=0$ est rationnelle. Or la condition $3)$ du théorème implique que toute solution de ($R_3$) est méromorphe, donc rationnelle, et par suite toute solution  $Y=-C^{-1}W'W^{-1}$ de l'équation $(R_1)$ est rationnelle, d'où le Théorème $7$.\\[0.2cm]
\underline{\bf Théorème $8$} :\\[0.2cm]
\textit{Si en tout point singulier $z\in Z$, l'opérateur $P$ est singulier régulier et les racines $k_1, \, ... \, , k_{2N}$ de son polynôme caractéristique $l(k)$ sont rationnelles et vérifient $k_i-k_j \notin \mathbb{Z}$ pour $i\neq j$, alors toute solution  de l'équation $(R_1)$ est algébrique.}\\[0.2cm]
\underline{Preuve}:\\[0.2cm]
On cherche une base de solutions de l'équation $(R_3)$ au voisinage de $x=a\in Z$ sous la forme de séries formelles 
$$W=x^k\sum\limits_{n\geq 0} W_ny^n,$$ avec $y=x-a$. \\
Elles doivent vérifier l'équation
$$\begin{array}{lll}
0=&Pu=\sum\limits_{n\geq 0}P(W_ny^{k+n})=L(k)W_0y^{k-n}+F_1(k)W_0y^{k-n+1}+F_2(k)W_0y^{k-n+2}+...\\
 &+L(k+1)W_1y^{k-n+1}+F_1(k+1)W_1y^{k-n+2}+...\\
\end{array}$$
avec $F_i(k)\in \mathbb{C}[k] \quad (i=0,1,2,...)$.\\
Si $k_i$ est un racine du polynôme $l(k)$, par hypothèse $l(k_i+j)\ne 0$ pour $j \in \mathbb{N}$, et donc on peut déterminer de manière unique tous les coefficients $W_n$ par récurrence et ainsi on obtient une solution $W(k_i)$. Cette solution $W(k_i)$ est convergente (Théorème $6$) et les $W(k_i)$ ($1\leq i\leq 2N$) forment une base de l'espace vectoriel des solutions de l'équation $Pu=0$ au voisinage du point $a$, donc la solution générale $Y=-C^{-1}(x)W'W^{-1}$ de ($R_1$) est algébrique, i.e les coefficients de cette matrices sont algébriques sur $\mathbb{C}(x)$.
Ceci termine la preuve du théorème.
\\[2.5cm]
\hspace{-0.6cm}\underline{\bf BIBLIOGRAPHIE}:\\[0.3cm]
[1] K.  Betina, \textit{Sur l'indice des opérateurs différentiels ordinaires}, Lecture Notes in Math.,\\  n° 1075, p.1-48 Springer-Verlag. \\[0.2cm]
[2] K. Betina, \textit{Indices polynômiaux et solutions polynômiales de certains systèmes de Pfaff singuliers}, Analysis 12 , 195-216. \\[0.2cm]
[3] K. Betina, \textit{Sur l'équation de Riccati}, Analysis 19, 29-50.\\[0.2cm]
[4] J.P. Bezivin, P. Robba, \textit{A new $p$-adic method for proving irrationality and transcendence results}, Annals of Maths, 129 (1989), 151-160 \\[0.2cm]
[5] H. Komatsu, \textit{On the index of differential operators}, Journal. Fac  Sci. Tokyo, IA (1971), 379-398.\\[0.2cm]
[6] B. Malgrange, \textit{Sur les points singuliers des équations différentielles}, L'enseignement Mathématique, t. XX (1974), 1-2.\\[2.5cm]
\end{document}